\documentclass[11pt]{article}
\usepackage{amssymb}
\usepackage{amsmath}

\newtheorem{theo}{Theorem}[section]

\newtheorem{prop}[theo]{Proposition}
\newtheorem{lemm}[theo]{Lemma}
\newtheorem{coro}[theo]{Corollary}
\newtheorem{rema}[theo]{Remark}
\newtheorem{Defi}[theo]{Definition}

\newtheorem{conj}[theo]{Conjecture}

\newcommand{\cqfd}
{%
\mbox{}%
\nolinebreak%
\hfill%
\rule{2mm}{2mm}%
\medbreak%
\par%
}
\newfont{\gothic}{eufb10}
\date{\empty}
\begin{document}
\title{Coniveau $2$ complete intersections and effective cones }
\author{Claire Voisin\\CNRS, Institut de math{\'e}matiques de Jussieu and IH\'{E}S} \maketitle \setcounter{section}{-1}
\section{Introduction}
\setcounter{equation}{0}
The goal of this paper is first of all to propose a strategy to attack the generalized Hodge conjecture for coniveau $2$ complete intersections,
and secondly to state a conjecture concerning the  cones of effective cycle classes in intermediate dimensions.
Our main results show that the generalized Hodge conjecture for coniveau $2$ complete intersections would follow from a particular
case of this effectiveness conjecture.

A rational Hodge structure of weight $k$ is given by a $\mathbb{Q}$-vector space $L$ together with a
Hodge decomposition
$$L_{\mathbb{C}}=\bigoplus_{p+q=k}L^{p,q}$$
satisfying Hodge symmetry
$$\overline{L^{p,q}}=L^{q,p}.$$
The coniveau of such a Hodge structure is the smallest integer $c$ such that
$L^{k-c,c}\not=0$.

When the Hodge structure comes from geometry, the notion of coniveau is conjecturally related to
codimension by the generalized Grothendieck-Hodge conjecture.
Suppose $X$ is a smooth complex projective variety and
$L\subset H^k(X,\mathbb{Q})$ is a sub-Hodge structure of coniveau $c$.
\begin{conj} \label{GHC} (cf \cite{grothendieck}) There exists a closed algebraic subset $Z\subset X$
 of codimension $c$ such that $L$ vanishes under the restriction
map  $H^k(X,\mathbb{Q})\rightarrow H^k(U,\mathbb{Q})$, where $U:=X\setminus Z$.
\end{conj}
Notice that it is a non trivial fact that the kernel of this restriction map is
indeed  a sub-Hodge structure of coniveau $\geq c$. This needs some arguments from mixed Hodge theory (see \cite{de2},
\cite{grothendieck} or \cite{voisinbook}, II, 4.3.2).
\begin{rema}{\rm  When a sub-Hodge structure $L$ has coniveau $c$, we can consider the generalized Grothendieck-Hodge conjecture for $L$
and for any coniveau $c'\leq c$, that is, ask whether $L$ vanishes on the complementary set of a closed algebraic subset
of codimension $c'$. When this holds, we will say  that $L$ satisfies the generalized Grothendieck-Hodge conjecture for coniveau $c'$.
}
\end{rema}
Consider a complete intersection $X\subset \mathbb{P}^n$ of $r$ hypersurfaces of degree
$d_1\leq \ldots\leq d_r$. By Lefschetz hyperplane sections theorem, the only interesting
Hodge structure is the Hodge structure on $H^{n-r}(X,\mathbb{Q})$, and in fact on the primitive part of it (that is the orthogonal of the restriction
of $H^*(\mathbb{P}^n,\mathbb{Q})$ with respect to the intersection pairing). We will say that $X$ has coniveau $c$
if the Hodge structure on $H^{n-r}(X,\mathbb{Q})_{prim}$ has coniveau $c$.

The coniveau of a complete intersection can be computed using Griffiths residues and the comparison of
pole order and Hodge filtration (see \cite{griffiths} or \cite{voisinbook}, II, 6.1.2). The result is as follows:
\begin{theo}\label{theogrif}  $X$ has coniveau $\ge c$ if and only if
$$n\geq \sum_id_i+(c-1)d_r.
$$
\end{theo}
For $c=1$, the result is obvious, as $coniveau(X)\geq1$ is equivalent to
$H^{n-r,0}(X)=H^0(X,K_X)=0$, that is $X$ is a Fano complete intersection.
In this case, the generalized Hodge-Grothendieck conjecture is known to be true, using the
correspondence between $X$ and its Fano variety of lines of lines $F$. Denoting by
$$\begin{matrix}&P&\stackrel{q}{\rightarrow} &X\\&p\downarrow&&\\&F&&
   \end{matrix}
$$
the incidence correspondence, where $p$ is the tautological $\mathbb{P}^1$-bundle on $F$,
one can show (see for example \cite{shimada}) that taking a $n-r-2$-dimensional complete intersection
$F_{n-r-2}\subset F$ and restricting $P$ to it, the resulting morphism of Hodge structures
$$q'_*\circ{p'}^*: H^{n-r-2}(F_{n-r-2},
\mathbb{Q})\rightarrow H^{n
-1}(X,
\mathbb{Q})$$
is
surjective, where $P'=p^{-1}(F_{n-r-2})$ and
$p',\,q'$ are the restrictions of $p,\,q$ to $P'$.
It follows that $H^{n-1}(X,\mathbb{Q})$ vanishes on the complement of the (singular) hypersurface
$q(P')\subset X$.

In the case of coniveau $2$, the numerical condition given by theorem \ref{theogrif} becomes
\begin{eqnarray}\label{estimate}n\geq \sum_id_i+d_r.
\end{eqnarray}
The geometric meaning of this bound is not so obvious. Furthermore the generalized Grothendieck-Hodge conjecture for coniveau
$2$
 is not known to hold in this range. For a general hypersurface $X$ of degree $d$ in $\mathbb{P}^n$, it is known
  to hold (for coniveau $2$) only when the dimension $n$ becomes much larger than the degree $d$, so that $X$
becomes covered by a family of planes (see \cite{evl}), (a slightly weaker condition has been obtained by Otwinowska \cite{otw}, see section \ref{applisec}).
The numerical range in which the result is known looks like
$$n\geq \frac{d^2}{4}
+O(d)$$ which is very different from the bound
$n\geq 2d$
of theorem \ref{theogrif}. For specific hypersurfaces of coniveau $\geq c$ defined by an equation
$F(X_0,\ldots,X_{n})$ of the form
$$F=F_1(X_0,\ldots,X_{d+s})+F_2(X_{d+s+1},\ldots,X_{2d+s})+\ldots + F_r(X_{(c-1)d+s+1}\ldots,X_{cd+s}),$$
(so $n=cd+s$),
 the generalized Hodge conjecture \ref{GHC} is proved in \cite{voisincras}
which proves more generally  that they satisfy the corresponding Bloch-Beilinson conjecture on Chow groups of hypersurfaces:
\begin{conj} \label{BB} If $coniveau(X)\geq c$, the Chow groups $CH_i(X)_{hom,\mathbb{Q}}$ of cycles homologous to zero modulo rational equivalence are trivial for $i\leq c-1$.
\end{conj}
Notice that the previously mentioned papers  \cite{evl}, \cite{otw}  are also devoted to the study of cycles of low dimension, and that the
generalized Grothendieck-Hodge conjecture is deduced from vanishings for these.
We refer to \cite{bloch-srinivas}, \cite{schoen} or \cite{voisinbook}, II, 10.3.2, for the proof that the Bloch-Beilinson conjecture implies
the generalized Hodge-Grothendieck conjecture.

We propose in this paper a strategy to prove the generalized Hodge conjecture for coniveau $2$ complete intersections, which does involve
the study of Chow groups and the construction of $2$-cycles (replacing the lines used in the coniveau $1$ case).

Our method is based on the following result which allows to give in section
\ref{intersec} a geometric proof of the numerical estimate (\ref{estimate}) for the coniveau $2$ property :
First of all, let us make the following definition:
\begin{Defi} \label{verymoving} A smooth $k$-dimensional subvariety $V\subset Y$,  where $Y$ is smooth projective, is {\rm very moving} if it  has the following property: through a generic point $y\in Y$, and given a generic
vector subspace $W\subset T_{Y,y}$ of rank $k$, there is a  deformation
$Y'$ of $Y$  which is smooth and passes through   $y$ with tangent space equal to $W$ at $y$.
\end{Defi}
 Let $X\subset \mathbb{P}^n$ be a generic complete intersection of multidegree $d_1\leq\ldots\leq d_r$.
For a generic section $G\in H^0(X,\mathcal{O}_X(n-\sum_id_i-1))$, consider
the subvariety $F_G\subset F$ of the variety of lines in $X_G$. Then by genericity, $F$ and $F_G$ are smooth of respective dimensions
$2n-2-\sum_id_i-r$ and $n-r-2$. One can show  that the deformations of $F_G$ are given by deformations of $G$.
\begin{theo} When $$n\geq \sum_id_i+d_r$$
the subvariety $F_G\subset F$ is ``very moving''.
\end{theo}

Cones of effective cycles have been very much studied in codimension $1$ or in dimension $1$
(cf  \cite{bdPP}), but
essentially nothing is known in intermediate (co)dimensions. Let us
say that an algebraic cohomology class is big if it belongs
to the interior of the effective cone.  In \cite{pet}, it is shown that when $dim\,W=1$, and
$W\subset V$  is moving and has ample normal bundle, its class $[W]$ is big. We will give here an example working in any dimension
$\geq 4$
and in codimension $2$, showing that in higher dimensions, a moving variety $W\subset V$ with ample normal bundle has not a necessarily
big class. Here by {\it moving}, we mean that a generic deformation of $W$ in $V$ may be imposed to pass
through a generic point of $V$.

We make the following conjecture for ``very moving'' subvarieties.

\begin{conj}\label{conjjbig} Let $V$ be smooth and  projective and let $W\subset V$ be a very moving
subvariety. Then the class $[W]$ of $W$ is big.
\end{conj}
This conjecture, in the case of codimension $2$  complete intersections in projective space, predicts rather mysterious effectiveness statements concerning projective bundles on Grassmannians.
 We will show this effectiveness result exactly in the same range appearing in \cite{otw}, which together with
 theorem \ref{main} below gives another proof of the fact that the generalized Hodge-Grothendieck
 for coniveau $2$ is satisfied in this numerical range.

We finally prove the following result (theorem \ref{mmain}) in  section \ref{applisec} :
\begin{theo}\label{main} Assume $n\geq \sum_id_i+d_r$ and the subvariety $F_G\subset F$ introduced above has a big class
(that is satisfies conjecture \ref{conjjbig}).
Then the complete intersections of multidegree $d_1\leq\ldots\leq d_r$ in $\mathbb{P}^n$ satisfy the generalized Hodge conjecture for coniveau $2$.
\end{theo}

Our method reproves the known results  concerning the generalized  Hodge conjecture
for coniveau $2$, that is proves it in the same range as \cite{otw}, but the spirit is very different, and the two methods lead in fact to a different statement, which we explain to conclude this introduction.
It is generally  believed that to solve the generalized Hodge conjecture for coniveau $c$ for $H^k(X,\mathbb{Q}),\,k=dim\,X$, one should produce a family of cycles
$(Z_{b})_{b\in B}, \,dim\,B=k-2c$,
of dimension $c$ in $X$, such that
the incidence family
$$\begin{matrix}&\mathcal{Z}&\stackrel{q}{\rightarrow} &X\\
&p\downarrow&&\\
&B&
\end{matrix}
$$
induces a surjective map
$$q_*p^*:H^{k-2c}(B,\mathbb{Q})\rightarrow H^k(X,\mathbb{Q}).$$
This obviously implies that $H^k(X,\mathbb{Q})$ vanishes away from $q(\mathcal{Z})$ and thus that the generalized  Hodge conjecture for coniveau $c$ is satisfied.
In general such cycles are provided by the proof that Chow groups of dimension $<c$ are small (cf \cite{bloch-srinivas} or \cite{voisinbook}, II, proof of theorem 10.31).

However  if the generalized Hodge conjecture holds for coniveau $c$ and for $H^k(X,\mathbb{Q})$, it does not  imply
the existence of such family, unless we also have a Lefschetz type conjecture satisfied. To see this more precisely, suppose the Hodge conjecture
holds true for $H^k(X,\mathbb{Q})$ and for  coniveau $c$. Then there exists a closed algebraic subset
$Z\subset X$ of codimension $c$ such that $H^k(X,\mathbb{Q})$ vanishes on $X\setminus Z$.
Introduce a desingularization
$\tau:\widetilde{Z}\rightarrow X$ of $Z$. Then the vanishing of $H^k(X,\mathbb{Q})$ on $X\setminus Z$ implies by strictness
of morphisms of mixed Hodge structures (see \cite{de2} or \cite{voisinbook}, II, 4.3.2) that
$$\tau_*:H^{k-2c}(\widetilde{Z},\mathbb{Q})\rightarrow H^k(X,\mathbb{Q})$$
is surjective.
Observe now that $dim\,\widetilde{Z}=k-c$. If the Lefschetz standard conjecture is satisfied by $\widetilde{Z}$, there exists a variety $B$ of dimension
$k-2c$ and a cycle $\mathcal{T}\subset B\times\widetilde{Z}$ of codimension $k-2c$ (a family of cycles on $\widetilde{Z}$ of dimension $c$ parameterized by
$B$), such that the map
$$q_*p^*:H^{k-2c}(B,\mathbb{Q})\rightarrow H^{k-2c}(\widetilde{Z},\mathbb{Q})$$
hence also the map
$$\tau_*\circ q_*p^*:H^{k-2c}(B,\mathbb{Q})\rightarrow H^k(X,\mathbb{Q})$$
are surjective, where $p$ and $q$ are the maps from $\mathcal{T}$ to $B$ and $\widetilde{Z}$ respectively.

Hence the parametrization of $H^k(X,\mathbb{Q})$ by algebraic cycles of dimension $c$ does not follow from the generalized Hodge conjecture for coniveau $c$, but also needs a
Lefschetz standard conjecture applied to a certain subvariety of $X$.

\section{A geometric interpretation of the coniveau $2$ condition \label{intersec}}
In this section, we will give a geometric interpretation of the
numerical condition
(\ref{estimate}), relating it to a  positivity property of a certain cycle class
on the variety of lines of the considered complete intersection.
We will also show how to deduce theorem \ref{theogrif} for coniveau $2$ from this positivity property.

Let thus $X$ be a generic complete intersection of multidegree $d_1\leq\ldots\leq d_r$ in
$\mathbb{P}^n$. Thus the variety of lines of $X$ is smooth
of dimension $2n-2-\sum_id_i-r$.
For $G$ a generic polynomial of degree $n-\sum_id_i-1$, let $X_G\subset X$ be the hypersurface
defined by $G$ and
$F_G\subset F$ the variety of lines contained in $X_G$. Thus $F_G$ is smooth of dimension
$$2n-2-\sum_id_i-r-(n-\sum_id_i)=n-r-2=dim\,X-2.$$
Recall the incidence diagram
$$\begin{matrix}&q:&P&\rightarrow X\\
&&\downarrow&\\
&&F&
\end{matrix}
$$
which induces for $n-\sum_id_i\geq0$ an injective morphism of Hodge structures:
$$p_*q^*:H^{n-r}(X,\mathbb{Q})\rightarrow H^{n-r-2}(F,\mathbb{Q}).$$
\begin{lemm}\label{elem} For a primitive cohomology class $a\in H^{n-l}(X,\mathbb{Q})_{prim}$,
the class $\eta:=p_*q^*a\in H^{n-r-2}(F,\mathbb{Q})$ satisfies the following
two properties (property \ref{1} will be used only in section \ref{applisec}):

\begin{enumerate}
 \item \label{1} (see  \cite{shimada}) $\eta$ is primitive with respect to the
Pl\"ucker polarization $l:=c_1(\mathcal{L})$ on $F$.
\item\label{2} $\eta$ vanishes on subvarieties $F_G$:

\begin{eqnarray}\label{van}
\eta_{\mid F_G}=0.
\end{eqnarray}

\end{enumerate}
\end{lemm}
{\bf Proof.} For the proof of the first statement, recall first that primitive cohomology
$H^{n-r-2}(F,\mathbb{Q})_{prim}$ is defined as the kernel of
$$\cup l^{n-\sum_id_i+1}:H^{n-r-2}(F,\mathbb{Q})\rightarrow H^{3n-2\sum_id_i-r}(F,\mathbb{Q})$$
because $dim\,F=2n-2-\sum_id_i-r$.
On the basis $U$ parameterizing smooth complete intersections $X$ such that $F$ is smooth of the right dimension,
the composed maps
$$\cup l^{n-\sum_id_i+1}\circ p_*q^*:H^{n-r}(X)_{prim}\rightarrow H^{n-r-2}(F)$$
give a morphism of local system.
The point is now the following: suppose $X$ degenerate to a generic  $X_0$ with one ordinary  double point $x_0$.
Then the family $Z_0$ of lines in $X_0$ passing through $x_0$ has dimension
$n-\sum_id_i$. It follows that it does not meet the generic intersection $K$ of $n-\sum_id_i+1$ members of the Pl\"ucker linear system
$\mid\mathcal{L}\mid$. Choose an $X_\epsilon$ which is generic and close enough to $X_0$. Then $X_\epsilon$ contains a
vanishing sphere $S_\epsilon$ which is arbitrarily close to $x_0$. Thus the cycle
$p_\epsilon(q_\epsilon^{-1}(S_\epsilon))\subset F_\epsilon$ does not meet a small perturbation
$K_\epsilon\subset F_\epsilon$ of $K$. It follows that if $\delta_\epsilon$ is the class of $S_\epsilon$ (well defined up to a sign depending on
a choice of orientation),
$\gamma_\epsilon:=p_{\epsilon*}q_\epsilon^*(\delta_\epsilon) \in H^{n-r-2}(F_\epsilon)_{prim}$ is primitive. Hence we conclude that $\delta_\epsilon$ belongs to the kernel of this
morphism of local systems. As the monodromy along $U$, acting on $H^{n-r}(X,\mathbb{Q})_{prim} $, acts transitively on
vanishing cycles (transported from the  $\delta_\epsilon$'s), and the later generate  $H^{n-r}(X,\mathbb{Q})_{prim} $ (cf \cite{voisinbook}, II, 3.2.2), it follows that
the morphism
$\cup l^{n-\sum_id_i+1}\circ p_*q^*:H^{n-r}(X)_{prim}\rightarrow H^{n-r-2}(F)$ is zero.

The second statement is elementary. Indeed, as $X_G\subset X$ is a smooth member of
$\mid\mathcal{O}_X(n-\sum_id_i-1)\mid$, by Lefschetz theorem, the restriction map
$$H^{n-r}(X,\mathbb{Q})_{prim}\rightarrow H^{n-r}(X_G,\mathbb{Q})$$
is zero.
Now we have the following commutative diagram:
$$\begin{matrix}&H^{n-r}(X,\mathbb{Q})_{prim}&\rightarrow &H^{n-r}(X_G,\mathbb{Q})&\\
&p_*q^*\downarrow&&p_{G*}q_G^* \downarrow&\\
&H^{n-r-2}(F,\mathbb{Q})&\rightarrow& H^{n-r-2}(F_G,\mathbb{Q})&
\end{matrix}
$$
where
$$\begin{matrix} &P_G&\stackrel{q_G}{\rightarrow}
& X_G\\
&p_G\downarrow&\\
&F_G&
\end{matrix}
$$ is the incidence diagram for $X_G$ and the horizontal maps are restriction maps.
Thus we have
$\eta_{\mid F_G}=p_{G*}q_G^*(a_{\mid X_G})=0$.
\cqfd
The main result of this section is the following:
\begin{theo} If $X$ is as above and  $n\geq \sum_id_i+d_r$, the subvarieties $F_G\subset F$, where $deg\,G=n-\sum_id_i-1$ are very moving
(cf definition \ref{verymoving}).
\label{propmain}
\end{theo}
Before giving the proof, let us use it, combined with lemma \ref{elem}, to give a geometric proof of  the numerical estimate (\ref{estimate}) of theorem \ref{theogrif} for coniveau $2$.
\begin{coro} (cf theorem \ref{theogrif}) If the inequality $n\geq \sum_id_i+d_r$ is satisfied, $X$ has coniveau $\geq 2$.
\end{coro}
{\bf Proof.} The fact that $X$ has coniveau $\geq2$ is equivalent (as we know already that $H^{n-r,0}(X)=0$) to the vanishing
$H^{n-r-1,1}(X)_{prim}=0$.
Let $a\in H^{n-r-1,1}(X)_{prim}$ and consider
$$\eta:=p_*q^*a\in H^{n-r-2,0}(F)=H^0(F,\Omega_F^{n-r-2}).$$
 Then $\eta=0$ iff $a=0$.
We use now statement \ref{2} in lemma \ref{elem}. This gives us
$$\eta_{\mid F_G}=0\,\,{\rm in}\,\,H^0(F_G,\Omega_{F_G}^{n-r-2})$$
for generic $G$.
On the other hand, theorem \ref{propmain} tells that for generic $G$,
the $r-2$ dimensional subvariety $F_G$ passes through a generic point $\Delta\in F$ with a generic tangent space.
It follows immediately that
$\eta_{\mid F_G}=0$ implies $\eta=0$.

\cqfd
{\bf Proof of theorem \ref{propmain}.}
We fix a line $\Delta\in F$. We want to study the differential of
the map $\phi$ which to
a polynomial $G$ vanishing on $\Delta$ associates the tangent space at $\Delta$ of $F_G$, assuming it has the right dimension.
What we will prove is the fact that the differential of $\phi$ is generically surjective when the bound is realized.
Let $W:=\phi(G)$, that is
$W\subset T_{F,\Delta}$ is the tangent space to $F_G$.
Then the differential of $\phi$ is a linear map
\begin{eqnarray}\label{dphig}d\phi(G):H^0(X,\mathcal{I}_{\Delta}(n-\sum_id_i-1))\rightarrow Hom\,(W,T_{F,\Delta}/W)
\end{eqnarray}
where the right hand side is the tangent space to the Grassmannian of rank $n-r-2$ dimensional subspaces
of $T_{F,\Delta}$.
We remark that  $d\phi(G)$ factors through
$$H^0(\Delta,\mathcal{I}_{\Delta,X}/\mathcal{I}_{\Delta,X}^2(n-\sum_id_i-1))=H^0(\Delta,N_{\Delta/X}^*(n-\sum_id_i-1)).$$
Furthermore, we observe that the natural map
$$H^0(X,\mathcal{I}_{\Delta,X}(n-\sum_id_i-1))\rightarrow H^0(\Delta,\mathcal{I}_{\Delta,X}/\mathcal{I}_{\Delta,X}^2(n-\sum_id_i-1))$$
is surjective. Indeed, the restriction map
$$H^0(\mathbb{P}^n,\mathcal{I}_{\Delta,\mathbb{P}^n}(n-\sum_id_i-1))\rightarrow
H^0(\Delta,\mathcal{I}_{\Delta,\mathbb{P}^n}/\mathcal{I}_{\Delta,\mathbb{P}^n}^2(n-\sum_id_i-1))$$
$$=H^0(\Delta,N_{\Delta/\mathbb{P}^n}^*(n-\sum_id_i-1))$$
is surjective. Furthermore the
the conormal exact sequence:
$$0\rightarrow\bigoplus_j\mathcal{O}_\Delta(-d_j)\rightarrow N_{\Delta/\mathbb{P}^n}^*\rightarrow N_{\Delta/X}^*\rightarrow 0$$
shows that the cokernel of the  map
$$H^0(\Delta, N_{\Delta/\mathbb{P}^n}^*(n-\sum_id_i-1))\rightarrow H^0(\Delta,N_{\Delta/X}^*(n-\sum_id_i-1))$$
is \begin{eqnarray}\label{obssp}H^1(\Delta,\bigoplus_j\mathcal{O}_\Delta(n-\sum_id_i-1-d_j)).
\end{eqnarray}

As $d_r=Sup\,\{d_j\}$, the  vanishing of the space (\ref{obssp})  is clearly implied (and in fact equivalent to) by  the inequality
$n-\sum_id_i-1-d_r\geq -1$ of (\ref{estimate}).

We now show the surjectivity of the map induced by $d\phi(G)$:
$$\overline{d\phi(G)}:H^0(\Delta,N_{\Delta/X}^*(n-\sum_id_i-1))\rightarrow Hom\,(W,T_{F,\Delta}/W).$$
Let $G\in H^0(X,\mathcal{I}_{\Delta/X}(n-\sum_id_i-1))$ be generic and $W=T_{F_G,\Delta}$.
Then the quotient $T_{F,\Delta}/W$ identifies to $H^0(\Delta,\mathcal{O}_\Delta(n-\sum_id_i-1))$. Indeed, $F_G$ is defined as the zero set of a
transverse section of $S^{n-\sum_id_i-1}\mathcal{E}$, and thus the normal bundle of $F_G$ in $F$ identifies to
$S^{n-\sum_id_i-1}\mathcal{E}_{\mid F_G}$, with fiber $H^0(\Delta,\mathcal{O}_\Delta(n-\sum_id_i-1))$ at $\Delta$.
Furthermore we  observe that as $\Delta$ is generic and $n-\sum_id_i>1$, the normal bundle
$N_{\Delta/X}$ is generated by sections. Of course $N_{\Delta/X}^*(n-\sum_id_i-1)$ is generated by global sections too.
 Finally, we note that by definition,
the space $W$ is the space $H^0(\Delta,N_{\Delta/X_G})$. As
$deg\,G=n-\sum_id_i-1$, $X_G$ is a generic Fano complete intersection of index $2$, and
it thus follows that  the normal bundle $N_{\Delta/X_G}$ is generically a direct sum of copies of $\mathcal{O}_\Delta$.

We have the following lemma:
\begin{lemm}\label{lemulti}1) The space $W=T_{F_G,\Delta}\subset H^0(\Delta,N_{\Delta/X})$ identifies to the kernel of the contraction map with
$$G\in H^0(\Delta,\mathcal{I}_{\Delta,X}/\mathcal{I}_{\Delta,X}^2(n-\sum_id_i-1))=H^0(\Delta,N_{\Delta/X}^*(n-\sum_id_i-1)),$$
with value in $H^0(\Delta,N_{\Delta/X}^*(n-\sum_id_i-1))$.

2) Using the inclusion
$W\subset H^0(\Delta,N_{\Delta/X})$, the map  $\overline{d\phi(G)}$
is induced by the contraction map between
$H^0(\Delta,N_{\Delta/X})$ and $H^0(\Delta,N_{\Delta/X}^*(n-\sum_id_i-1))$, with value in $H^0(\Delta,\mathcal{O}_\Delta(n-\sum_id_i-1))$.
\end{lemm}
Postponing the proof of this lemma, we now conclude as follows.
Let $G\in H^0(\Delta,N_{\Delta/X}^*(n-\sum_id_i-1))$ be generic and let
$W=Ker \,d\phi(G)$, where $d\phi(G)$ has been identified in lemma \ref{lemulti} to contraction by
$G$, with value in $H^0(\Delta,\mathcal{O}_\Delta(n-\sum_id_i-1))$. We have to show that
the map given by contraction
$$H^0(\Delta,N_{\Delta/X}^*(n-\sum_id_i-1))\rightarrow Hom\,(W,H^0(\Delta,\mathcal{O}_\Delta(n-\sum_id_i-1)))$$
is surjective. This problem concerns now vector bundles on $\Delta=\mathbb{P}^1$:
we have a vector bundle $E$ of rank $s=n-r-1$ and degree $k=n-1-\sum_id_i$ on $\Delta=\mathbb{P}^1$,  such that
$E^*(k)$ is generated by global sections.
We choose a generic element $G$ of $H^0(\Delta,E^*(k))$. We know that
$G$ gives
a surjective map
$E\rightarrow \mathcal{O}_\Delta(k)$ with kernel $K$ which is a trivial vector bundle. Denote by $W\subset H^0(\Delta,E)$ the kernel
of the contraction map with $G$, with value in $H^0(\mathcal{O}_\Delta(k))$; thus $W=H^0(\Delta,K)$ and we have to show that the
contraction map induces a surjective map:
$$H^0(\Delta,E^*(k))\rightarrow Hom\,(W,H^0(\Delta,\mathcal{O}_\Delta(k)).$$

We consider now
the composed map
$$H^0(\Delta,E^*(k))\rightarrow H^0(\Delta,K^*(k))\rightarrow Hom\,(H^0(\Delta,K),H^0(\Delta,\mathcal{O}_\Delta(k)))$$
and want to show that it is surjective. The first map is surjective as its cokernel is $H^1(\Delta,\mathcal{O}_\Delta)=0$.
The second map is surjective exactly when $H^0(\Delta,K(-1))=0$ which follows from the fact that $K$ is trivial.

\cqfd
{\bf Proof of lemma \ref{lemulti}.} 1)  We have the inclusions
$$\Delta\subset X_G\subset X,$$
which give the normal bundles exact sequence:
$$0\rightarrow N_{\Delta/X_G}\rightarrow N_{\Delta/X}\stackrel{dG}{\rightarrow} \mathcal{O}_\Delta(n-\sum_id_i-1)\rightarrow0.$$
The space $W$ is by definition the space $H^0(\Delta,N_{\Delta/X_G})$ hence by the exact sequence above, it is the kernel of
$$dG:H^0(\Delta,N_{\Delta/X})\rightarrow H^0(\Delta,\mathcal{O}_\Delta(n-\sum_id_i-1)).$$
We just have to observe that the map $dG$ above is nothing but contraction with the image of $G$ in $\mathcal{I}_\Delta/\mathcal{I}_\Delta^2(n-\sum_id_i-1)=N_{\Delta/X}^*(n-\sum_id_i-1)$, which follows from the construction of the normal bundles sequence.

\vspace{0.5cm}

2) It is an immediate consequence of 1) and the following more general statement: Let
$\phi:V\rightarrow H$ be a surjective map and let $W:=Ker\,\phi$.
A small  deformation of $\phi$, given by an element of $Hom\,(V,H)$ gives a deformation of $W\in Grass(w,V)$. Thus we have a natural rational map
$$a: Hom\,(V,H)\rightarrow Grass(w,V)$$
where $w=rk\,W$.
 Then the differential
$h\in Hom\,(V,H)\mapsto da(h)\in Hom\,(W,V/W)=Hom\,(W,H)$ of this map is simply
the map
$$h\in Hom\,(V,H)\mapsto h_{\mid W}.$$
This last statement follows from the standard construction of the isomorphism $T_{Grass(w,V)}\cong Hom\,(W,V/W)$. This concludes the proof of 2).
\cqfd
\section{A conjecture on cones of effective cycles}
Let $Y$ be a smooth projective complex variety and let $Alg^{2k}(Y)\subset H^{2k}(Y,\mathbb{R})$ be the vector subspace
generated by classes of codimension $k$ algebraic cycles. This vector space contains
the effective cone $$E^{2k}(Y)\subset Alg^{2k}(Y)$$
generated by classes of effective cycles.

\begin{Defi}   A class $\alpha\in Alg^{2k}(Y)$ is said to be big
if $\alpha$ is an interior point of $E^{2k}(Y)$.
\end{Defi}
If $h=c_1(H)$ where $H$ is an ample line bundle on $Y$, $h^k$ belongs to the interior of the cone
$E^{2k}(Y)$. Indeed, for any effective cycle $Z\subset Y$ of codimension $k$,
the class $N h^k-[Z]$ is effective for $N$ large enough. Applying this to a basis of $Alg^{2k}(Y)$ consisting of effective cycles,
one concludes that for some open set $U\subset Alg^{2k}(Y)$ containing $0$, $h^k-U\subset E^{2k}(Y)$.

Thus we get the following:
\begin{lemm}\label{lemmepsilon} A class $\alpha\in Alg^{2k}(Y)$ is big if and only if, for some
$\epsilon >0$, $\alpha-\epsilon h^k\in E^{2k}(Y)$.
\end{lemm}
In the divisor case, the effective cone, or rather its closure, the pseudo-effective cone, is now well understood
by  work of Boucksom, Demailly, Paun and Peternell \cite{bdPP}. The case of higher codimension is not understood at all. To start with,
there is the following elementary result for divisors, which we will show to be wrong in codimension $2$ and higher.
\begin{lemm} \label{trivialdiv} Let $D$ be an effective divisor on $Y$, and assume that
$\mathcal{O}_Y(D)_{\mid D}$ is ample. Then $[D]$ is big.
\end{lemm}
{\bf Proof.}
Indeed, one shows   for any divisor $E$ on $Y$,  by vanishing on $D$, that the natural map
$$H^1(Y,\mathcal{O}_Y((n-1)D-E))\rightarrow H^1(Y,\mathcal{O}_Y(nD-E))$$ is
surjective for large $n$, hence must be an isomorphism for large $n$. Thus the restriction map
$$H^0(Y, \mathcal{O}_Y(nD-E))\rightarrow H^0(D,\mathcal{O}_D(nD-E))$$ is surjective for
large $n$. The space $H^0(D,\mathcal{O}_D(nD-E))$ is non zero
for $n$ large, by ampleness of $\mathcal{O}_D(D)$, and it thus follows that the space $H^0(Y,\mathcal{O}_Y(nD-E))$ is non zero
too. Taking for $E$ an ample divisor and applying lemma \ref{lemmepsilon} gives the result.
\cqfd
Note that when $D$ is smooth, $\mathcal{O}_D(D)$ is also the normal bundle $N_{D/Y}$.
Let us now construct   examples of codimension $2$ smooth subvarieties $W\subset Y$ with ample normal bundle,
such that $[W]$ is not big.

We start from a generic  hypersurface $Z$ of degree $d$ in $\mathbb{P}^{d+2}$, with $d\geq3$ and denote by $Y$ its variety of lines.
Let now $\mathbb{P}^{d+1}\subset \mathbb{P}^{d+2}$ be an hyperplane, and let
$W\subset Y$ be the subvarieties of lines contained in $Z':=Z\cap \mathbb{P}^{d+1}$. By genericity, $W$ and $Y$ are
smooth, and $W$ is the zero set of a transverse section of the tautological rank $2$ bundle
$\mathcal{E}$ on $F$, with fiber $H^0(\Delta,\mathcal{O}_\Delta(1))$ at the point $\Delta\in F$. Thus the normal bundle
of $W$ in $Y$ is isomorphic to $\mathcal{E}_{\mid W}$.
This vector bundle is globally generated on $W$, and it is ample if the
natural map
$$\mathbb{P}(\mathcal{E})\rightarrow Z'$$
is finite to one. This is indeed the case for generic $Z'$ by \cite{jiang}.
Thus $W$ has ample normal bundle.
Let us now show:
\begin{prop} \label{contreex} The class $[W]\in Alg^4(Y,\mathbb{R})$ is not big.
\end{prop}
{\bf Proof.} The $d+1$-dimensional hypersurface $Z$ has $H^{d,1}(Z)_{prim}\not=0$, and this provides using the incidence diagram
(cf section \ref{intersec}) a non zero holomorphic form $\eta\in H^{d-1,0}(Y)$ vanishing on $W$.
Note that $dim\,W=d-1$.
Assume that $[W]$ is big. By lemma \ref{lemmepsilon}, one has
\begin{eqnarray}\label{egalite}[W]=\epsilon h^2+e,
\end{eqnarray}
where $e=\sum_i e_i[E_i]$ is a class of an effective codimension $2$ cycle with real coefficients and $h$ is the class of an ample divisor.
Let us integrate the form $i^{d-1}(-1)^{\frac{(d-1)(d-2)}{2}}\eta\wedge\overline{\eta}$ on both sides. As $\eta_{\mid W}=0$, the left hand side is zero.
On the other hand, we have by the second Hodge-Riemann relations (for holomorphic forms):
$$\int_{E_i}i^{d-1}(-1)^{\frac{(d-1)(d-2)}{2}}\eta\wedge\overline{\eta}\geq 0,\,\int_Yh^2\cup i^{d-1}(-1)^{\frac{(d-1)(d-2)}{2}}\eta\wedge\overline{\eta}>0.$$
Thus the integral on the right is $>0$, which is a contradiction.
\cqfd

Note that the variety $W$ also has the property that it is moving, that is, its deformations sweep out $Y$. The proof of lemma
\ref{contreex} shows however that it is not very moving (see section \ref{intersec}). Let us make the following conjecture :
\begin{conj} \label{conjbig} Let $W\subset Y$ be a smooth variety which is very moving.
then $[W]$ is big.
\end{conj}
\begin{rema} {\rm We could weaken slightly the conjecture above by considering varieties which are very moving and have ample normal bundle.  Indeed, the varieties we will consider in the sequel and for which we would need
to know that their class is big not only are very moving (cf theorem \ref{propmain}) but also have ample normal bundle.}
\end{rema}
The conjecture is true for divisors by lemma \ref{trivialdiv}.
We also have:
\begin{prop}  Conjecture \ref{conjbig} is true for curves.
\end{prop}
{\bf Proof.}
 Indeed, by Mo\"\i shezon-Nakai criterion,
the dual cone to the effective cone of curves is the ample cone. A point
$c=[C]$  which is in the boundary of the effective cone of curves
must thus have zero intersection with a nef class $e$, that is a point in the boundary of the ample cone.
Consider now a one-dimensional family $(C_t)_{t\in B}$ of curves deforming $C$ and passing through a given
generic point $y\in Y$. Such a family exists as $C$ is very moving.
Let $\tau: \Sigma\rightarrow X$ be a desingularization of this family. Then $\Sigma$ contains two divisor classes
$c$ and $\sigma$ which correspond to the fibers of the map  $\Sigma\rightarrow B$ and to the section
of this map given by the point $y$. The nef class
 $\tau^*e\in H^2(\Sigma,\mathbb{R})$ satisfies
  $$<\tau^*e,c>=0,<\tau^*e,\tau^*e>\geq 0$$
  while the class $ c$ satisfies $<c,c> =0$, where the $<,>$ denotes the intersection product on $H^2(\Sigma,\mathbb{R})$. The Hodge index theorem then tells that
  the classes $c$ and $\tau^* e$ must be proportional. But as $<c,\sigma> >0$ and $<\tau^* e,\sigma>=0$, it follows that in fact
  $\tau^*e=0$.
The contradiction comes from the following: as $e\not=0$ there is a curve $\Gamma\subset X$, which we may assume to be
in general position, such that $deg_\Gamma e>0$. Now, because $C$ is very moving, there exist for each $\gamma\in \Gamma$
curves $C_\gamma$ which are deformations of $C$ and pass through
$y$ and $\gamma$. We can thus choose our one dimensional family
to be parameterized by a cover $r:\Gamma'\rightarrow \Gamma$, in such a way that
for any $\gamma'\in \Gamma'$, $r(\gamma')\in C_{\gamma'}$.
But then the surface $\tau(\Sigma)$ contains the section $\gamma'\mapsto\gamma'\in C_{r(\gamma')}$ which sends
 via $\tau$ onto $\Gamma$, and as $deg_\Gamma e>0$, we conclude that
$\tau^*e\not=0$.
\cqfd
Another example where conjecture \ref{conjbig} holds is the following:
\begin{lemm} \label{schubert}Let $V\subset G(1,n)$ be a very moving subvariety.
Then $[V]$ is big.
\end{lemm}
{\bf Proof.} The cone of effective cycles on $G(1,n)$ is very simple because (cf \cite{manivel})
we have for each
pairs of complementary dimensions dual basis of the cohomology which are generated by classes of Schubert cycles, which are effective, and furthermore,
as the tangent space of the Grassmannian is globally generated,
 two effective cycles of complementary dimension on $G(1,n)$ have non negative intersection (cf \cite{fula}).
 If $z$ is an effective class on $G(1,n)$, write $z=\sum_i\alpha_i\sigma_i$ where $\sigma_i$ are
 Schubert classes of dimension equal to $dim\,z$. Then $\alpha_i=z\cdot\sigma_i^*\geq0$.
 Thus the cone of effective cycles is the cone generated by classes of Schubert cycles.
 It follows from this argument that for the class $[V]$ to be big, it suffices that
 $$V\cdot \sigma_i^*>0$$ for every Schubert cycle $\sigma_i^*$ of complementary dimension.
 These  inequalities  are  now implied by the fact that $V$ is very moving. Indeed,
 choose a smooth point $x\in \sigma_i^*$. Then we can choose a deformation $V'$ of $V$ passing through $x$,
 smooth at $x$ and meeting transversally $\sigma_i^*$ at $x$.
 Thus there is a non zero contribution to $V\cdot\sigma_i^*$ coming from the point $x$.
 The intersection $V\cap\sigma_i^*$ could be non proper away from $x$, but because
 the tangent space of the Grassmannian is generated by global sections, each
 component of the intersection $V\cap\sigma_i^*$ has a non negative contribution to $V\cdot\sigma_i^*$ by \cite{fula}.
 Thus $V\cdot\sigma_i^*>0$.
\cqfd
\subsection{On the conjecture \ref{conjbig} for the varieties $F_G$}
We finally turn to the  study of conjecture \ref{conjbig} for the subvarieties $F_G\subset F$. For $X$ a generic  complete intersection
of multidegree $d_1\leq\ldots\leq d_r$ in $\mathbb{P}^n$,
we have seen that $F_G\subset F$ is very moving, where $deg\,G=n-\sum_id_i-1$.
We have now the following:
\begin{theo}\label{partitheo} The class $[F_G]\in H^{2n-2\sum_id_i}(F,\mathbb{Q})$ is big
when
\begin{eqnarray}\label{bound}3n-4-\sum_i\frac{(d_i+1)(d_i+2)}{2}\geq n-r-2.
\end{eqnarray}
\end{theo}
Before proving this, let us explain the geometric meaning of this bound. The number
$3n-6-\sum_i\frac{(d_i+1)(d_i+2)}{2}$ is simply the expected dimension of the family of planes contained in $X$.
The number $3n-4-\sum_i\frac{(d_i+1)(d_i+2)}{2}$ is thus the expected dimension of the family $Z'$ of lines contained in a plane contained in $X$.
On the other hand, the dimension $n-r-2$ appearing on the right is the dimension of the varieties $F_G$.

\vspace{0.5cm}

{\bf Proof of theorem \ref{partitheo}} First of all, by Lemma \ref{debutrec} proved below, it suffices to consider
the case where
$$3n-4-\sum_i\frac{(d_i+1)(d_i+2)}{2}= n-r-2.$$
What we shall do is to compute the class of the subvariety $Z'$ of
$F$ introduced above, which we
describe as follows: let $G_2\rightarrow G$ be the partial flag manifold parameterizing pairs $(Q,\Delta),\,\Delta\subset Q$, where $\Delta$ is a line in $\mathbb{P}^n$ and
$Q$ is a plan in $\mathbb{P}^n$. Let $\pi:F_2\rightarrow F$ be the subvariety of $G_2$ consisting of such pairs
with $\Delta\in F$. Thus $\pi:F_2\rightarrow F$ is a projective bundle of relative dimension $n-2$.
 We consider the variety $Z\subset F_2$ consisting of pairs $(\Delta,P)\in F_2$ where $P$ is also  contained in $X$. It has dimension
 $3n-4-\sum_i\frac{(d_i+1)(d_i+2)}{2}$ which by assumption is equal to $ n-r-2$.  Thus
 $Z'=\pi(Z)$  has dimension $n-r-2$  and its   codimension
 is $n-\sum_id_i$.

 Theorem \ref{partitheo} is an immediate consequence of the following two
 statements, (lemma \ref{leok} and proposition \ref{anotherprop}).
 \cqfd

 \begin{lemm} \label{leok} Assume there exists a class $P\in H^{2n-2\sum_id_i}(F,\mathbb{Q})$,
 which has the following properties:
 \begin{enumerate}
 \item $P$ is an effective class on $F$.
 \item \label{item2} $P$ can be written as
 $$P=-\epsilon l^{n-\sum_id_i}+c_2R,$$
 where $R$ is any algebraic class of $F$.
 \end{enumerate}
 Then the class $c_{n-\sum_id_i}(S^{n-\sum_id_i-1}\mathcal{E})\in H^{2n-2\sum_id_i}(F,\mathbb{Q})$ is big.
 \end{lemm}
 \begin{prop} \label{anotherprop} The class $[Z']$ (which is by definition effective) of the variety $Z'$ defined above  is given by a polynomial expression $P=P(l,c_2)$ satisfying property \ref{item2} of lemma \ref{leok}.
 \end{prop}
 We used above the notation  $c_2=c_2(\mathcal{E})\in H^4(F,\mathbb{Q})$ and $l=c_1(\mathcal{E})\in H^2(F,\mathbb{Q})$, $\mathcal{E}$ being
 as before
 the (restriction to $F$ of) the tautological rank $2$ vector bundle with fiber $H^0(\Delta,\mathcal{O}_\Delta(1))$ at $\Delta\in F$.

\vspace{0.5cm}

{\bf Proof of lemma \ref{leok}} Introduce a formal splitting of $\mathcal{E}$ or equivalently formal roots  $x,y$ of its Chern polynomial, so that
 $$c_2=xy,\,l=x+y.$$
 Then $S^{n-\sum_id_i-1}\mathcal{E}$ has for formal roots the expressions
 $k x+(n-\sum_id_i-1-k)y$, with $0\leq k\leq n-\sum_id_i-1$.
 Thus we get
 $$c_{n-\sum_id_i}(S^{n-\sum_id_i-1}\mathcal{E})=\prod_{k=0}^{n-\sum_id_i-1}k x+(n-\sum_id_i-1-k)y$$
 which one can rewrite as
 $$c_{n-\sum_id_i}(S^{n-\sum_id_i-1}\mathcal{E})=(n-\sum_id_i)^2xy\prod_{k=1}^{n-\sum_id_i-2}k x+(n-\sum_id_i-1-k)y.$$
 We claim that the class $Q:=\prod_{k=1}^{n-\sum_id_i-2}k x+(n-\sum_id_i-1-k)y$ is big on the Grassmannian
 hence a fortiori on $F$.
To see this, we apply the argument of lemma \ref{schubert}. We just have to show that
$$Q(x,y)\cdot \sigma>0$$ for all Schubert cycles of   dimension  $n-\sum_id_i-2$ on $G(1,n)$.
We now use the fact that
$$(n-\sum_id_i-1)^2c_2Q=c_{n-\sum_id_i}(S^{n-\sum_id_i-1}\mathcal{E}).$$
Any Schubert cycle of  dimension $n-\sum_id_i-2$ on $G(1,n)$ is of the form
$\sigma=c_2\sigma'$ where $\sigma'$ is a Schubert cycle of dimension $n-\sum_id_i$ on $G(1,n)$.
One then has
$$(n-\sum_id_i-1)^2Q(x,y)\cdot \sigma=(n-\sum_id_i)^2c_2Q(x,y)\cdot\sigma'=c_{n-\sum_id_i}(S^{n-\sum_id_i-1}\mathcal{E})\cdot\sigma'.$$
One shows easily that for any Schubert cycle
$\sigma'$ of dimension $n-\sum_id_i$, one has
$$ c_{n-\sum_id_i}(S^{n-\sum_id_i-1}\mathcal{E})\cdot \sigma'>0$$
unless $\sigma'=\sigma_{0,n-\sum_id_i+1}$ is the Schubert cycle of lines passing through a point and contained
in a linear space of dimension $n-\sum_id_i+1$.
But in this case, $c_2\sigma'=0$. This proves the claim.

Now we proved that $c_{n-\sum_id_i}(S^{n-\sum_id_i-1}\mathcal{E})=c_2Q$, where $Q$ is big on $F$. Assume
that there is a class $P$ of codimension $n-\sum_id_i$ on $G$ which is effective
and of the form $-\epsilon l^{n-\sum_id_i}+c_2R$, with $R$ algebraic.
As the class $Q$ is big on $F$, it is in the interior of the effective cone of $F$,
and thus for some small $\epsilon'$
$$Q-\epsilon'R$$
is effective on $F$.
Thus $c_2(Q-\epsilon'R)$ is also effective on $F$ and we get:
$$c_{n-\sum_id_i}(S^{n-\sum_id_i-1}\mathcal{E})=c_2Q=\epsilon'c_2R+E,$$
with $E$ effective. Replacing $c_2R$ by $P+\epsilon l^{n-\sum_id_i}$, where $P$ is effective, we get:
$$c_{n-\sum_id_i}(S^{n-\sum_id_i-1}\mathcal{E})=E+\epsilon'P+\epsilon'\epsilon l^{n-\sum_id_i},$$ where $E+\epsilon'P$ is effective. By lemma \ref{lemmepsilon}, $c_{n-\sum_id_i}(S^{n-\sum_id_i-1}\mathcal{E})$ is then big.

\cqfd
{\bf Proof of proposition \ref{anotherprop}.} We give for simplicity the proof
for $r=1$ and thus denote $d_r=d$ with $n=2d$. Let $\pi:F_2\rightarrow F$ be as above. $F_2$ is a $\mathbb{P}^{n-2}$-bundle over $F$.
Let $P_2\rightarrow F_2$ be the universal plane parameterized by $F_2$. As $P_2$ is sent naturally to $\mathbb{P}^n$, this $\mathbb{P}^2$-bundle
admits a natural polarization $\mathcal{O}(1)$. Then
$$P_2=\mathbb{P}(\mathcal{F})$$
for a certain rank $3$ vector bundle on $F_2$.
As for each pair $(\Delta,Q)\in F_2$, one has $\Delta\subset Q$, one has a natural surjective
  map of bundles on $F_2$
$$\mathcal{F}\rightarrow \pi^*\mathcal{E}\rightarrow 0.$$
Let $\mathcal{H}$ be its kernel and let $h:=c_1(\mathcal{H})\in H^2(F_2,\mathbb{Q})$.
The $\mathbb{P}^{n-2}$-bundle $F_2$ over $F$ is polarized by the line bundle $\mathcal{L}_2:=det\,\mathcal{F}$ coming from the Grassmannian
$G(2,n)$ via the natural map $F_2\rightarrow G(2,n),\,(\Delta,Q)\mapsto Q$.
From the exact sequence
\begin{eqnarray}
\label{exactseq}0\rightarrow\mathcal{H}\rightarrow\mathcal{F}\rightarrow\pi^*\mathcal{E}\rightarrow0,
\end{eqnarray}
one deduces the exact sequence
$$0\rightarrow\mathcal{H}\otimes S^{d-1}\mathcal{F}\rightarrow S^d\mathcal{F}\rightarrow \pi^*S^d\mathcal{E}\rightarrow0.$$
Observe now that the defining equation $f$ of $X$ induces a section
of $S^d\mathcal{F}$ which by definition vanishes in $S^d\mathcal{E}$. Thus we get a section
$\tilde{f}$ of $\mathcal{H}\otimes S^{d-1}\mathcal{F}$ such that $Z=V(\tilde{f})$.
It follows (together with some transversality arguments implying that $\tilde{f}$ vanishes transversally) that
$$[Z]=c_N(\mathcal{H}\otimes S^{d-1}\mathcal{F}),\, N=rank \, S^{d-1}\mathcal{F}=\frac{d(d+1)}{2}.$$
It remains to compute
$$[Z']=\pi_*[Z].$$
It is clear that $[Z']$ is a polynomial in $l$ and $c_2$.
As we are interested only in showing that the coefficient of $l^{n-d}$ in this polynomial is negative, we can do formally as if
$\mathcal{E}=\mathcal{L}\oplus\mathcal{O}_F$.
As $\mathcal{H}\otimes S^{d-1}\mathcal{F}$ is filtered with successive quotients isomorphic to
$$\mathcal{H}^{\otimes i}\otimes S^{d-i}\mathcal{E},\,i=1\ldots,\,d$$
we get
$$[Z]=\prod_{i=1}{d}c_{d-i+1}(\mathcal{H}^{\otimes i}\otimes S^{d-i}\mathcal{E})$$
and modulo $c_2$ this is equal to
$$\prod_{i=1}^{d}\prod_{j=0}^{d-i}(ih+(d-i-j)l).$$
Observe now that by the exact sequence (\ref{exactseq}), the $\mathbb{P}^{n-2}$-bundle $P_2$   polarized by $\mathcal{H}=\mathcal{L}_2\otimes\pi^*\mathcal{L}^{-1}$ is isomorphic to
$\mathbb{P}(\mathcal{K}_1)$,
where $\mathcal{K}_1$ is the kernel of the evaluation map
$V\rightarrow \mathcal{E}$, so that
 we have the exact sequence
$$0\rightarrow \mathcal{K}_1\rightarrow V\otimes \mathcal{O}_F\rightarrow \mathcal{E}\rightarrow0.$$
By definition of Segre classes, we then have
$$\pi_*h^{n-2+i}=s_i(\mathcal{K}_1^*)=c_i(\mathcal{E}^*).$$
As we compute modulo $c_2$, we get:
$$\pi_*h^{n-2+i}=0,\,i\not=0,\,1$$
$$\pi_*h^{n-2}=1,\,p_*h^{n-1}=-l.$$

Hence it follows that we have the equality
$$[Z']=\pi_*(\prod_{i=1}^{d}\prod_{j=0}^{d-i}(ih+(d-i-j)l))=(\alpha_{n-2}-\alpha_{n-1})l^{n-d}\,\,mod.\,\,c_2,$$
where we write
$$M:=\prod_{i=1}^{d}\prod_{j=0}^{ d-i}(ih+(d-i-j)l))=\sum_{i\leq N}\alpha_ih^il^{N-i}.$$

The important point is now the following:
Let us factor in $M$ all the terms corresponding to
$j=d-i$. Then
$$M=d!h^d\prod_{i\geq1,j\geq 1,i+j\leq d}(ih+jl).$$
Let $M':=\prod_{i\geq1,j\geq 1,i+j\leq d}(ih+jl)=\sum_{i\leq N-d}\beta_ih^il^{N-d-i}$.
Then we have
$$\alpha_{n-2}-\alpha_{n-1}=d!(\beta_{n-d-2}-\beta_{n-d-1}).$$
Thus we have to show that $\beta_{n-d-2}-\beta_{n-d-1}<0$.
But now we observe that the degree of the homogeneous polynomial $M'$
is $N-d$ with $N=n-2+n-d$. Hence $deg\,M'=2n-2d-2$. The polynomial $M'$ is  symmetric
of degree $2n-2d-2$ in $l$ and $h$ and it suffices to show that its coefficients
$\beta_i$ are strictly increasing
in the range $i\leq n-d-1=\frac{deg\,M'}{2}$.
This can checked at hand or proved by geometry using the fact that
$$(d!)^2(xy)^dM'(x,y)=\prod_{1\leq i\leq d}c_{i+1}(S^i\mathcal{E}),$$
if $x,\,y$ are the formal roots of $\mathcal{E}$, and that the numbers
$(d!)^2(\beta_{i}-\beta_{i+1}),\,2i\leq 2n-2d-2$ can be interpreted as intersection numbers of $\prod_{1\leq i\leq d}c_{i+1}(S^i\mathcal{E})$
with a Schubert cycle.
\cqfd
Finally, let us rephrase  conjecture \ref{conjbig} in the case of subvarieties $F_G\subset F$, in terms of geometry of the Grassmannian
$G:=G(1,n)$ of lines in
$\mathbb{P}^n$.
For simplicity, let us consider the case of hypersurfaces of degree $d$ in $\mathbb{P}^n$. The bound (\ref{estimate}) then becomes
$$n\geq 2d.$$
Let $V:=H^0(\mathbb{P}^n,\mathcal{O}_{\mathbb{P}^n}(1))$ with symmetric power  $S^dV=H^0(\mathbb{P}^n,\mathcal{O}_{\mathbb{P}^n}(d))$. Recall that we denote by
$\mathcal{E}$ the rank $2$ vector bundle with fiber $H^0(\Delta,\mathcal{O}_\Delta(1))$ on $G(1,n)$. There is a natural evaluation map
\begin{eqnarray}\label{evaln}S^dV\otimes \mathcal{O}_G\rightarrow S^d\mathcal{E}\end{eqnarray}
with kernel we denote by $\mathcal{K}_d$. The projective bundle
$$\pi:\mathbb{P}(\mathcal{K}_d^*)\rightarrow G(1,n)$$
(where we use the Grothendieck notation) thus parameterizes pairs
$(f,\Delta)$ such that $f_{\mid \Delta}=0$.
The natural projection $\rho:\mathbb{P}(\mathcal{K}_d^*)\rightarrow\mathbb{P}(S^dV^*)$, which to $(\Delta,f)$ associates $f$, has for fibre over $f$ the
variety $F_f$ of lines contained in the hypersurface $X_f$ defined by $f$.

Now we consider the subvarieties $F_{f,G}\subset F_f$, where $G$ has degree $n-d-1$. These are defined as zero sets of generic, hence transverse, sections of
$S^{n-d-1}\mathcal{E}$. Thus their cohomology class
is given by
$$[F_G]=c_{n-d}(S^{n-d-1}\mathcal{E}_{\mid F}).$$
Summarizing and applying lemma \ref{lemmepsilon}, conjecture \ref{conjbig} tells that for some $\epsilon>0$, (a multiple of) the class
$\pi^*(c_{n-d}(S^{n-d-1}\mathcal{E})-\epsilon l^{n-d})$  is effective on the generic fiber of the map
$\rho:\mathbb{P}(\mathcal{K}_d^*)\rightarrow\mathbb{P}(S^dV^*)$. An elementary Hilbert scheme argument then shows that we can put this in family,
and combining this with the description of the cohomology ring of
a projective bundle, this gives us the following reformulation of conjecture \ref{conjbig} in this case: we  denote by
$H$ the class $c_1(\mathcal{O}(1))$ on the projective bundle $\mathbb{P}(\mathcal{K}_d^*)$.
\begin{conj}\label{subbigconj} Assume $n\geq 2d$. Then for some algebraic class $$\alpha'\in H^{2n-2d-2}(\mathbb{P}(\mathcal{K}_d^*),\mathbb{Q}),$$ and for some small
$\epsilon>0$, the class
$$\pi^*(c_{n-d}(S^{n-d-1}\mathcal{E})-\epsilon l^{n-d})+H\alpha'$$
is effective on $\mathbb{P}(\mathcal{K}_d^*)$. In other words, the class $c_{n-d}(S^{n-d-1}\mathcal{E})$ belongs to the interior
of the cone of degree $2n-2d$ classes $\alpha$ on $G$ such that $\pi^*\alpha+H\alpha'$ is effective on $\mathbb{P}(\mathcal{K}_d^*)$ for some
algebraic class $\alpha'\in H^{2n-2d-2}(\mathbb{P}(\mathcal{K}_d^*),\mathbb{Q})$.
\end{conj}
Theorem \ref{partitheo} proves conjecture \ref{subbigconj} in the range
$3n-4-\frac{(d+1)(d+2)}{2}\geq n-3$, that is
$$n\geq \frac{1}{2}(\frac{(d+1)(d+2)}{2}+1).$$
\begin{rema}{\rm More generally, suppose we have a variety $Y$ and a vector bundle $\mathcal{F}\rightarrow Y$, with associated projective bundle
$\pi:\mathbb{P}(\mathcal{F})\rightarrow Y$. Let $h=c_1(\mathcal{O}_{\mathbb{P}(\mathcal{F})}(1))$. We can then introduce
the convex cone $E^{2k}(Y,\mathcal{F})$ consisting of classes $c\in Alg^{2k}(Y)$ such that $\pi^*c+hc'\in E^{2k}(\mathbb{P}(\mathcal{F}))$,
for some $c' \in Alg^{2k-2}(\mathbb{P}(\mathcal{F}))$. There is an obvious inclusion
$$E^{2k}(Y)\subset E^{2k}(Y,\mathcal{F}).$$
When $\mathcal{F}$ is ample, the whole of $Alg^{2k}(Y)$ is contained in $E^{2k}(Y,\mathcal{F})$ and it follows that
$E^{2k}(Y)$ is contained in the interior of $E^{2k}(Y,\mathcal{F})$. At the opposite, if $\mathcal{F}$ is trivial, it is obvious that
$$E^{2k}(Y)= E^{2k}(Y,\mathcal{F}).$$
In our case, the class
$c_{n-d}(S^{n-d-1}\mathcal{E})$ is effective and even numerically effective but not in the interior of the effective cone of $Alg^{2n-2d}(G(1,n))$. The bundle
$\mathcal{K}_d^*$ is generated by sections but not ample.

}
\end{rema}
\begin{rema}{\rm The conjecture cannot be improved by replacing $c_{n-d}(S^{n-d-1}\mathcal{E})$ by
$c_{n-d-1}(S^{n-d-2}\mathcal{E})$. Indeed,  on a generic Fano variety of lines  $F$ in a degree $d$ hypersurface
$X$, the class
$[F_G]=c_{n-d-1}(S^{n-d-2}\mathcal{E})$ of the subvariety of lines in $X_G,\,deg\,G=d-n-2$ is not big. In fact,  $dim\,F_G=n-2$ and $F_G$ has $0$ intersection with
the variety $F_x$ of lines through a generic point $x\in X$, which has dimension $n-d-1=dim\,F-n+2$.
As $F_x$ is moving, this implies that $[F_G]$ is not big. The same argument works for any class in $H^{2n-2d-2}(F,\mathbb{Q})$ which is divisible by
$c_2$.}
\end{rema}
  The crucial case for conjecture \ref{subbigconj} is the case where $n=2d$. Indeed, we have the following:
  \begin{prop}\label{debutrec} If the conjecture \ref{subbigconj} is true for a given pair $(n,d)$, $n\geq2d$, then it is true for the pair
  $(n+1,d)$.
\end{prop}
{\bf Proof.} Let us denote by $\mathcal{K}_d^{n+1}$ the kernel of the evaluation map
(\ref{evaln}) on $G(1,n+1)$, for  $V=H^0(\mathbb{P}^{n+1},\mathcal{O}_{\mathbb{P}^{n+1}}(1))$.
We assume that the conjecture is satisfied for $n, \,d$ with $n\geq 2d$ and show that it is satisfied
for $n+1,d$.

Choose a pencil of hyperplanes $(H_t)_{t\in\mathbb{P}^1}$ on $\mathbb{P}^{n+1}$, with base-locus $B$.
Then $G(1,n+1)$ contains the Pl\"ucker hyperplane section $G_B$ consisting of lines meeting $B$.
it  is singular and admits a natural desingularization
$$G'\rightarrow G_B\subset G(1,n+1)$$
where $G'=\sqcup _{t\in \mathbb{P}^1}G(1,H_t)$.
We pull-back  to
$G'$ the projective bundle $\mathbb{P}({\mathcal{K}_d^{n+1}}^*)$, and denote the result by
$$P':=G'\times _{G(1,n+1)} \mathbb{P}({\mathcal{K}_d^{n+1}}^*)$$
$$j:P'\rightarrow \mathbb{P}({\mathcal{K}_d^{n+1}}^*),$$
We denote by  $\pi_{n+1}:\mathbb{P}({\mathcal{K}_d^{n+1}}^*)\rightarrow G(1,n+1)$  the structural map,
and let $\pi'_{n+1}:=\pi_{n+1}\circ j:P'\rightarrow G(1,n+1)$. We have the following properties:
\begin{enumerate}\item\label{1sept5} $j_*[P']=\pi_{n+1}^*l\in H^2(\mathbb{P}({\mathcal{K}_d^{n+1}}^*,\mathbb{Q})$.
\item\label{2sept5} If $n':P'\rightarrow \mathbb{P}^1$ is the natural map, and
$k=c_1(\mathcal{O}_{\mathbb{P}^1}(1))$,
$$j_*({n'}^*k)=\pi_{n+1}^*c_2\in H^4(\mathbb{P}({\mathcal{K}_d^{n+1}}^*),\mathbb{Q}).$$
\end{enumerate}
Observe also that the cohomology of $P'$ is generated by $j^*H^*(\mathbb{P}({\mathcal{K}_d^{n+1}}^*),\mathbb{Q})$ and
${n'}^*k$, because for  each fiber $P'_t$ of $n'$, the restriction map
$$j_t^*:H^*(\mathbb{P}({\mathcal{K}_d^{n+1}}^*),\mathbb{Q})\rightarrow H^*(P'_t,\mathbb{Q})$$
is surjective.

Introduce $$P'':=\sqcup_{t\in \mathbb{P}^1}\mathbb{P}({\mathcal{K}_d^{n,t}}^*)$$
where $\mathbb{P}({\mathcal{K}_d^{n,t}}^*)$ is the variety $\mathbb{P}({\mathcal{K}_d^{n}}^*)$ for the hyperplane $H_t$.
There is a natural rational linear projection of projective bundles :
$$\phi:P'\dashrightarrow P''$$
which to $(\Delta,t,X)$ associates $(\Delta,t,X\cap H_t)$.  Denote by
$$\pi''_{n+1}:P''\rightarrow G(1,n+1),\,n'':P''\rightarrow \mathbb{P}^1$$
the natural maps. Then we have
$$\pi''_{n+1}\circ \phi=\pi'_{n+1}=\pi_{n+1}\circ j,\,n'=n''\circ\phi.$$

Let $h$ be the class $c_1(\mathcal{O}_{\mathbb{P}({\mathcal{K}_d^{n+1}}^*)}(1))$.
The projective bundle $P''\rightarrow G''$ has a natural polarization $h''$ such that
$$\phi^*h''=j^*h.$$
 The induction assumption is that for each
fiber $P''_t$ of $n''$,
a class of the form ${\pi''_{n+1}}^*(-\epsilon l^{n-d}+c_{n-d}(S^{n-d-1}\mathcal{E}))+h'' c'$ with $\epsilon>0$ is effective on $P''_t$. Putting this in family,
we conclude that
on $P''$, a class of the form
$${\pi''_{n+1}}^*(-\epsilon l^{n-d}+c_{n-d}(S^{n-d-1}\mathcal{E}))+h'' c'+{n''}^*k c''$$
is effective.
The fact that the indeterminacy locus of $\phi$ has high codimension then shows that
$$\phi^*({\pi''}_{n+1}^*(-\epsilon l^{n-d}+c_{n-d}(S^{n-d-1}\mathcal{E}))+h'' c'+{n''}^*k c'')$$
\begin{eqnarray}\label{derder}=j^*({\pi^*_{n+1}}(-\epsilon l^{n-d}+c_{n-d}(S^{n-d-1}\mathcal{E})))+j^*h \phi^*c'+{n'}^*k\phi^*c''
\end{eqnarray}
is also effective on $P'$.
As noted above, the class $\phi^*c''$ comes from a class on $\mathbb{P}({\mathcal{K}_d^{n+1}}^*)$,
$$\phi^*c''=j^*c'''.$$
Applying $j_*$ to (\ref{derder}), the properties \ref{1sept5} and \ref{2sept5} above, and the projection formula, we conclude that
a class of the form
$$l({\pi^*_{n+1}}(-\epsilon l^{n-d}+c_{n-d}(S^{n-d-1}\mathcal{E}))+h c')+({\pi^*_{n+1}}c_2)c'''$$
$$={\pi_{n+1}^*}(-\epsilon l^{n-d+1}+lc_{n-d}(S^{n-d-1}\mathcal{E}))+h lc'+({\pi^*_{n+1}}c_2)c'''$$
is effective on $ \mathbb{P}({\mathcal{K}_d^{n+1}}^*)$.

As $c_{n-d}(S^{n-d-1}\mathcal{E}))$ is divisible by $c_2$, we conclude
  that the assumptions of lemma \ref{leok} are  satisfied
  on the fibers of
 $$\rho:\mathbb{P}({\mathcal{K}_d^{n+1}}^*)\rightarrow \mathbb{P}(S^dV^*),$$
 where $h$ vanishes. Thus  by lemma \ref{leok}, conjecture \ref{subbigconj} is satisfied
 for $n+1,d$.

\cqfd

\section{On the generalized Hodge conjecture for coniveau $2$ complete intersections \label{applisec}}
We prove now the following result which motivated our interest in theorem \ref{propmain} and conjecture \ref{conjbig}.
 Let $X$ be a generic complete intersection of multidegree $d_1\leq\ldots\leq d_r$ in $\mathbb{P}^n$, and assume the bound
(\ref{estimate}) holds, that is
$$n\geq \sum_id_i+d_r.
$$
By theorem \ref{propmain}, the varieties $F_G\subset F$ are very moving.
\begin{theo} \label{mmain} If the varieties $F_G$ satisfy conjecture \ref{conjbig}, that is $[F_G]$ is big,
then the generalized Hodge conjecture for coniveau $2$ is satisfied by $X$.
\end{theo}
{\bf Proof.} Assume that $[F_G]$ is big. Then it follows by lemma \ref{lemmepsilon} that for some positive large integer
 $N$  and for some effective cycle $E$ of codimension $n-\sum_id_i$ on $F$, one has:
\begin{eqnarray}\label{encoreentier}N[F_G]= l^{n-\sum_id_i}+[E].
\end{eqnarray}
(Here we could work as well with real coefficients, but as we want actually to do geometry on $E$, it is better if $E$ is a true cycle.)
Now we recall from lemma \ref{elem} that
for $a\in H^{n-r}(X)_{prim}$, $\eta=p_*q^*a\in H^{n-r-2}(F)$ is primitive with respect to $l$ and furthermore vanishes on $F_G$, with
$dim\,F_G=n-r-2$.
Let us assume that
$a\in H^{p,q}(X)_{prim}$ and integrate $(-1)^{\frac{k(k-1)}{2}}i^{p-q}\eta\cup\overline{\eta},\,k=p+q-2=n-r-2$ over
both sides in (\ref{encoreentier}). We thus get
$$0=\int_F(-1)^{\frac{k(k-1)}{2}}i^{p-q}l^{n-\sum_id_i}\cup\eta\cup\overline{\eta}+\int_E(-1)^{\frac{k(k-1)}{2}}i^{p-q}\eta\cup\overline{\eta}.$$
As $\eta$ is primitive, and non zero if $a$ is non zero, by the second Hodge-Riemann bilinear relations (cf \cite{voisinbook}, I, 6.3.2),
we have $ \int_F(-1)^{\frac{k(k-1)}{2}}i^{p-q}l^{n-\sum_id_i}\cup\eta\cup\overline{\eta} >0$.
It thus follows that
$$\int_E(-1)^{\frac{k(k-1)}{2}}i^{p-q}\eta\cup\overline{\eta}<0.$$
Let $\widetilde{E}=\sqcup \widetilde{E_j}$ be a desingularization of the support of $E=\sum_jm_jE_j,\,m_j>0$.
Thus we have
$$\sum_jm_j\int_{\widetilde{E_j}}(-1)^{\frac{k(k-1)}{2}}i^{p-q}\eta\cup\overline{\eta}<0.$$
It thus follows that there exists one $E_j$ such that
\begin{eqnarray}\label{inegpasprim}\int_{\widetilde{E_j}}(-1)^{\frac{k(k-1)}{2}}i^{p-q}\eta\cup\overline{\eta}<0.
\end{eqnarray}
Choose an ample divisor   $H_j$ on each $\widetilde{E_j}$.
By the second Hodge-Riemann bilinear relations, inequality \ref{inegpasprim} implies that
${\eta}_{\mid \widetilde{E_j}}$ is not primitive with respect to the polarization given by $H_j$, that is
$\eta\cup [H_j]\not=0$ and in particular
$${\eta}_{\mid H_j}\not=0.$$

In conclusion, we proved that the composed map
$$H^{n-r}(X)_{prim}\stackrel{p_*q^*}{\rightarrow}H^{n-r-2}(F)\rightarrow \bigoplus H^{n-r-2}(H_j)$$
is injective, where the second map is given by restriction.
If we dualize this, recalling that $dim\,H_j=n-r-3$, we conclude that
$$\bigoplus H^{n-r-4}(H_j)\rightarrow H^{n-r}(X,\mathbb{Q})_{prim}$$
is surjective, where we consider the pull-backs of the incidence diagrams to $H_j$
$$\begin{matrix} &P_j&\stackrel{q_j}{\rightarrow}&X\\
&p_j\downarrow&&\\
&H_j&&
\end{matrix}
$$ and the map is the sum of the maps
$q_{j*}p_j^*$, followed by orthogonal projection onto primitive cohomology. As for $n-r$ even and $n-r\geq3$, the image of
$\sum _jq_{j*}p_j^*$ also contains the class $h^{\frac{n-r}{2}}$, (up to adding if necessary to the $H_j$ the class of a linear section of $F$)
it follows immediately that
the map
$$\sum _jq_{j*}p_j^*:\bigoplus H^{n-r-4}(H_j)\rightarrow H^{n-r}(X,\mathbb{Q})$$
is also surjective, if $n-r\geq3$ (the case $n-r=2$ is trivial).
This implies  that $H^{n-r}(X,\mathbb{Q})$ is supported on the $n-r-2$-dimensional variety
$\cup _jq_j(P_j)$, that is vanishes on $X\setminus \cup _jq_j(P_j)$. The result is proved.

\cqfd
This theorem combined with
theorem \ref{partitheo} allows us to reprove and generalize the previously known results concerning the generalized Hodge conjecture for coniveau $2$ hypersurfaces
(see \cite{otw}). Indeed the bound (\ref{bound}) of theorem \ref{partitheo} can also be reinterpreted as follows (in the case of hypersurfaces):
the generic hypersurface of degree $d$ in $\mathbb{P}^{n+1}$ is swept-out by plans.
In \cite{otw}, A. Otwinowska proves that this implies the triviality of $CH_1(X)_{\mathbb{Q},hom}$ for
$X$ any hypersurface of degree $d$ in $\mathbb{P}^{n}$ (note that it is likely
that the same can be done for complete intersections as well).
From Bloch-Srinivas argument (see \cite{bloch-srinivas}, \cite{voisinbook}, II, proof of theorem 10.31) one can also deduce from this that
$H^{n-r}(X)_{prim}$ vanishes on the complementary set of a closed algebraic subset of $X$ of codimension $2$, that is the generalized Hodge conjecture
for coniveau $2$ holds for $X$.

\end{document}